\documentclass{amsart}

\usepackage{amssymb}

\usepackage{graphicx}

\usepackage{upref}

\usepackage{hyperref}

\usepackage{eucal} 

\usepackage{enumerate} 

\newtheorem{theorem}{Theorem}[section]
\newtheorem{lemma}[theorem]{Lemma}
\newtheorem{prop}[theorem]{Proposition}
\newtheorem{cor}[theorem]{Corollary}

\theoremstyle{definition}
\newtheorem{definition}[theorem]{Definition}
\newtheorem{example}[theorem]{Example}

\theoremstyle{remark}
\newtheorem{remark}[theorem]{Remark}

\numberwithin{equation}{section}

\newcommand{\set}[1]{\left\{#1\right\}}
\newcommand{\minus}{\smallsetminus}

\newcommand{\tree}{\mathsf{T}}
\newcommand{\F}{\mathsf{F}}

\newcommand{\child}{\textup{\textsf{c}}}
\newcommand{\parent}{\textup{\textsf{p}}}

\newcommand{\VP}{\mathsf{v}^{\parent}}
\newcommand{\VC}{\mathsf{v}^{\child}}
\newcommand{\vc}{\mathsf{v}^{\child}}
\newcommand{\vp}{\mathsf{v}^{\parent}}

\newcommand{\bC}{\mathbb{C}}
\newcommand{\bR}{\mathbb{R}}
\newcommand{\bN}{\mathbb{N}}
\newcommand{\bZ}{\mathbb{Z}}

\newcommand{\bD}{\mathbb{D}}
\newcommand{\bQ}{\mathbb{Q}}

\newcommand{\bdr}{\partial}
\newcommand{\dist}[2]{\text{dist}\,(#1,#2)}

\newcommand{\MRev}[2]{\href{http://www.ams.org/mathscinet-getitem?mr=MR{#1}}{MR~{#2}}}

\DeclareMathOperator{\Mod}{mod}

\newcommand{\Kjulia}{\CMcal{K}}
\newcommand{\julia}{\CMcal{J}}
\newcommand{\puz}{\CMcal{P}}

\begin{document}

\title[Return times of polynomials as meta-{F}ibonacci numbers]
{Return times of polynomials\\
as meta-{F}ibonacci numbers}

\author{Nathaniel D. Emerson}

\address{Department of Mathematics\\ University of Southern California\\ Los Angeles, California 90089}

\email{nemerson@usc.edu}



\date{September 30, 2008}

\begin{abstract}
We consider generalized closest return times of a complex polynomial of degree at least two. Most previous studies on this
subject have focused on the properties of polynomials with particular return times, especially the Fibonacci numbers. We study
the general form of these closest return times. The main result of this paper is that these closest return times are
meta-Fibonacci numbers.  In particular, this result applies to the return times of a principal nest of a polynomial. Furthermore,
we show that an analogous result holds in a tree with dynamics that is associated with a polynomial.
\end{abstract}

\maketitle

\section{Introduction}

Consider the dynamics of a polynomial $f: \bC \to \bC$ (see \cite{CG} for example).  The most intuitive sequence of \emph{closest
return times} of a point $z_0$ under iteration by $f$ is defined as follows: Let $n_1 = 1$ and define $n_{k+1}$ as the least
integer such that $ |f^{n_{k+1}}(z_0) - {z_0} | < |f^{n_{k}}(z_0) - {z_0}|$. M.~Lyubich and J.~Milnor proved that there exists a
real quadratic polynomial such that the closest return times of its critical point are the Fibonacci numbers
\cite{Lyubich-Milnor}. However, these closest return times are not generally invariant if we conjugate $f$ by an affine map.
Thus, it is natural to consider more general types of return times.

The return times associated with a principal nest of a polynomial are affine invariants. B.~Branner and J.~Hubbard studied the
dynamics of polynomials with a disconnected Julia set and exactly one critical point with bounded orbit.  They showed that the
Fibonacci numbers occur as the return times of a {principal nest} of certain of these polynomials \cite[Ex.\ 12.4]{BH92}.

The above-mentioned studies are typical of previous work on return times in complex dynamics. They considered a specific sequence
of return times, the Fibonacci numbers, and derived properties of polynomials with the specified return times.  In contrast, the
main theorem of this paper is a general result about the form of return times of a complex polynomial.  There are no previously
published results about the general form of closest return times of complex polynomials.

We study polynomial dynamics by associating a polynomial of degree at least $2$ with a \emph{tree with dynamics} (Definition
\ref{defn: Twd - abstract} and \cite{E03}). We introduce a \emph{return nest} (Definition \ref{defn: return nest}) of a
polynomial, which is a generalization of a \emph{principal nest} (Definition \ref{defn: principal nest} and
\cite{Lyubich-Quadratics}).   Return nests have a combinatorial analogue in a tree \linebreak

\newpage
\noindent with dynamics. Our main theorem is that the closest return times of any return nest are meta-Fibonacci numbers (Theorem
\ref{thm: Main}).  By \emph{meta-Fibonacci numbers} we mean a sequence given by a Fibonacci-type recursion, where the recursion
varies with the index (see \cite{CCT} for an overview).

Meta-Fibonacci numbers have not previously been considered in the context of complex dynamics.  We introduce them by recalling
some generalizations of Fibonacci numbers. The Fibonacci numbers are recursively defined by adding the previous two terms of the
sequence: $u_k = u_{k-1} + u_{k-2}$. Adding the previous three terms yields the \emph{Tribonacci numbers}. If we add the previous
$r$ terms, we obtain \emph{$r$-generalized Fibonacci numbers} (``\emph{$r$-bonacci numbers}'') \cite{Miles}. The meta-Fibonacci
numbers that we consider in this paper are defined by the following recursion. Let $r$ be a function of $k$, and add the previous
$r(k)$ terms:

\begin{equation}\label{eq: r(k) equation}
     n_k = \sum_{j=1}^{r(k)} n_{k-j}.
\end{equation}
Given $r: \bZ^+\to \bZ^+$ with $r(k) \leq k $ for all $k$, we choose an initial condition $n_0$, and recursively define $n_k$ by
Equation \ref{eq: r(k) equation} for $k \geq 1$. We call the resulting sequence $(n_k)_{k= 0}^{\infty}$ a \emph{variable-$r$
meta-Fibonacci sequence generated by $r(k)$} \cite[Def.\ 1.1]{E06}. We present examples in Section \ref{sect: meta-Fib seqs}. In
order to describe the return times of polynomials, we need to allow for the possibility that $r(k)$ is arbitrarily large.  Hence,
we define $r(k)$ and $n_k$ for all integers.

\begin{definition}\cite[Def.\ 5.1]{E06}  \label{defn: r(k)-bonacci}
Let $r: \bZ \to \bZ^+$ and let $(n_k)_{k \in \bZ}$ be a double sequence of real numbers. We say $(n_k)$ is an \emph{extended}
variable-$r$ meta-Fibonacci sequence generated by $r(k)$ if Equation \ref{eq: r(k) equation} holds for all $k \in \bZ$. For
brevity, we write ``\emph{$r(k)$-bonacci numbers}.''
\end{definition}

We describe the type of return times that we will consider. Let $f$ be a polynomial of degree at least 2.  We can form a puzzle
of $f$ (see \S\ref{sect: poly TwD}), which decomposes the complex plane into topological disks called \emph{puzzle pieces} of
$f$. We consider  dynamically defined subsequences of a sequence of nested puzzle pieces.  Let $(P_l)_{l \in \bZ}$ be a sequence
of nested of puzzle pieces of $f$. A \emph{return nest} is a sub-nest $(P_{l(k)})_{k \in \bZ}$ such that
\[
    f^{n_k}(P_{l(k)}) = P_{l(k-1)} \quad \text{for all $k \in \bZ$},
\]
where $n_k = \min \set{n \geq 1: f^n(P_{l(k)}) = P_m \text{ for some } m}$. We call $(n_k)_{k \in \bZ}$ the \emph{return times}
of the return nest.  In some cases, we need to modify this definition for one value of $k$ (Definition \ref{defn: return nest}).
A {principal nest} \cite[\S3.1]{Lyubich-Quadratics} is a special type of return nest. The following is our main theorem.

\begin{theorem}\label{thm: Main}
The sequence of return times of any return nest of a polynomial is an extended variable-$r$ meta-Fibonacci sequence.
\end{theorem}

Let us fix some notation. Let $f$ be a complex polynomial of degree at least $2$. We say a point is \emph{persistent} if it has
bounded orbit under $f$ and \emph{escapes} otherwise. The \emph{filled Julia set of $f$} is the set of all persistent points of
$f$, and we denote it by  $\Kjulia(f)$. The \emph{Julia set of $f$} is the topological boundary of $\Kjulia(f)$, and we denote it
by $\julia(f) $.  By a classical result of Fatou and Julia, the Julia set of $f$ is connected if and only if all critical points
of $f$ are persistent.

Our main tool for studying polynomial dynamics is the combinatorial system of a tree with dynamics (Definition \ref{defn: Twd -
abstract}).  A puzzle decomposes the plane. A tree with dynamics codes this decomposition and keeps track of the key features of
the dynamics in a simpler setting. A polynomial with a disconnected Julia set has a canonical Branner-Hubbard puzzle \cite{BH92}.
We use this decomposition to define a canonical tree with dynamics (Definition \ref{defn: TwD of Poly w/Disc} and \cite{E03}).
When $f$ is a polynomial with a connected Julia set, we decompose the plane into a {Yoccoz puzzle} \cite{Hubbard_Loc_Con}. We
introduce the tree with dynamics associated with a Yoccoz puzzle (Definition \ref{defn: TwD of Poly w/Con}).  Thus, trees with
dynamics are a single combinatorial system for studying the dynamics of polynomials with either connected or disconnected Julia
sets. Theorem \ref{thm: Main} follows from an analogous result for trees with dynamics (Theorem \ref{thm: ret times of ret chains
Fib}).

Puzzles have been used in complex dynamics for some time.  Branner and Hubbard defined a canonical dynamic decomposition of the
plane for polynomials with exactly one escaping critical point \cite[Ch.\ 1.1]{BH92}.  Later, Branner clarified the structure of
the Branner-Hubbard puzzle \cite{Branner}.  The author generalized this decomposition, and defined the tree with dynamics for any
polynomial with a disconnected Julia set \cite[\S3]{E03}. J.~C.~Yoccoz introduced a technique to decompose the the plane using
the dynamics of a polynomial with a connected Julia set, which is now called a ``Yoccoz puzzle''. For a quadratic polynomial, the
construction was first published by Hubbard \cite[$\S5$]{Hubbard_Loc_Con}. The general case was described by J. Kiwi
\cite[\S12]{Kiwi_Rat_Rays}).

Trees with dynamics were introduced by R.~P{\'e}rez-Marco \cite{PM_DCS}. The first substantive results using trees with dynamics
to study polynomials with a disconnected Julia sets were obtained by the author \cite{E_Thesis, E03}.  In the case of a
polynomial with a disconnected Julia set, Theorem \ref{thm: Main} is part of the author's thesis \cite{E_Thesis}.  The connected
case has not previously been published.

The remainder of this paper is organized as follows.  In Section \ref{sect: TwD}, we define trees with dynamics abstractly, and
prove some preliminary results about them.  We prove that a combinatorial version of Theorem \ref{thm: Main} holds in any tree
with dynamics (Theorem \ref{thm: ret times of ret chains Fib}).   We then give a construction of a combinatorial analogue of a
return nest.  Finally, we describe the initial conditions satisfied by the return times of any return nest.  In Section
\ref{sect: meta-Fib seqs}, we study the properties of $r(k)$-bonacci sequences and present some examples.   We derive estimates
on the growth of $n_k$ in terms of $r(k)$.  In Section \ref{sect: poly TwD}, we prove our main theorem and define a tree with
dynamics for a polynomial. We consider abstract puzzles with compatible dynamics, which are generalized Markov partitions. We
show there is a tree with dynamics associated with each abstract puzzle. This definition establishes a correspondence between a
puzzles and trees with dynamics.  Theorem \ref{thm: Main} follows from this correspondence. We outline the construction of a
puzzle for a polynomial with either a disconnected Julia set or with a connected Julia set. We finish by noting some properties
of a tree with dynamics of any polynomial.

\section{Trees with Dynamics}\label{sect: TwD}

In this section, we study trees with dynamics and prove a combinatorial version of our main theorem.  The proof is
straightforward once the necessary machinery is in place.  We define a general tree with dynamics (\S\ref{subs: TwD - prelim}). A
return chain is an analogue of a return nest in a tree with dynamics (Definition \ref{defs: ret chain}). The dynamics of return
nests are the subject of \S\ref{subs: Twd - Dyn of Ends}, where we prove several important lemmas. We prove a version of Theorem
\ref{thm: Main} for return chains (Theorem \ref{thm: ret times of ret chains Fib}) in \S\ref{subs: TwD - Main}.  We give
necessary and sufficient conditions for a tree with dynamics to have a return chain (\S\ref{subs: TwD - Ret Chain}). Finally, we
consider a class of trees with dynamics that includes all trees with dynamics of polynomials (\S\ref{subs: Twd - Rooted Trees}).
The return times of these trees with dynamics satisfy certain restrictions (Proposition \ref{prop: n_k = 1 for k < K}).

\subsection{Preliminaries}\label{subs: TwD - prelim}

A \emph{tree} is a countable connected graph with every circuit trivial.   We say two vertices of a graph are \emph{adjacent} if
there is an edge between them. We only consider trees with a particular type of order on their vertices.

\begin{definition} \label{defn: tree}
A \emph{genealogical tree} is a tree $\tree$ such that each vertex $\mathsf{v} \in \tree$ is associated with a unique adjacent
vertex $\vp$, the \emph{parent} of $\mathsf{v}$.  Every vertex adjacent to $\mathsf{v}$, except $\vp$, is called a \emph{child}
of $\mathsf{v}$ and denoted by $ \vc$.
\end{definition}

\begin{figure}[hbt]
\begin{center}

\includegraphics{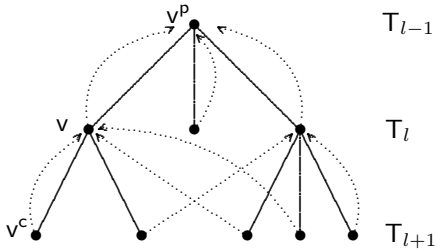}

\caption{A tree with dynamics with $H  = 1$.}
\label{fig: TwD}
\end{center}
\end{figure}

In this paper, by ``tree'' we mean genealogical tree.  We use the symbol $\tree$ to represent both the tree and its vertex set;
the edge set is left implicit. We use \textsf{sans serif} symbols for trees and objects associated with trees.  Our convention in
drawing trees is that a parent is above its children (see Fig.\ \ref{fig: TwD}). So $\VP$ is above $\mathsf{v}$ and any $\VC$ is
below $\mathsf{v}$.  When it is necessary to distinguish between children of $\mathsf{v}$ we use the notation
$\mathsf{v}^{\child_i}$. We say $\mathsf{v}$ is an \emph{ancestor} of $\mathsf{v}'$ if there are vertices $\mathsf{v}_0, \dots,
\mathsf{v}_n$ such that $\mathsf{v} = \mathsf{v}_0$, $\mathsf{v}' = \mathsf{v}_n$, and $\mathsf{v}_{i-1} =
\mathsf{v}_i^{\parent}$ for $i = 1, \dots, n$.  We say $\mathsf{v}''$ is a \emph{descendant} of $\mathsf{v}$ if $\mathsf{v}$ is
an ancestor of $\mathsf{v}''$.

Definition \ref{defn: tree} also defines a genealogical \emph{graph}.  A genealogical graph is a tree if and only if any two
vertices have a common ancestor and there is no vertex $\mathsf{v}$ such that $\vp$ is a descendant of $\mathsf{v}$.

Let $\tree$ be a tree such that $\tree = \bigcup_{l \in \bZ} \tree_l$, where $\tree_l = \set{\mathsf{v} \in \tree: \vp \in
\tree_{l-1}}$ for each $l$.  We call each $\tree_l$ a \emph{level} of $\tree$. We can inductively partition any tree $\tree$ into
levels, and this partition is unique up to a choice of $\tree_0$.  Thus, the levels of $\tree$ are well-defined up to indexing.
Hereafter we will assume that any tree has its levels indexed in some manner.

We consider all infinite paths in the tree that move from parent to child.  Let $\bN$ denote the non-negative integers.

\begin{definition}

Let $\tree$ be a tree.  An \emph{end} of ${\tree}$ is a sequence $\boldsymbol{\mathsf{x}}=(\mathsf{x}_l)_{l \in \bN}$, where
$\mathsf{x}_l \in \tree_l$ and $\mathsf{x}_{l-1} = \mathsf{x}_l^{\parent}$ for all $l$. An \emph{extended end} is the analogous
double sequence  $\boldsymbol{\mathsf{x}}=(\mathsf{x}_l)_{l \in \bZ}$.
\end{definition}

A natural metric for the extended ends of $\tree$ is a Gromov metric:
\[
    \dist{\boldsymbol{\mathsf{x}}}{\boldsymbol{\mathsf{y}}} = \gamma^{-L}, \quad
    L = \max \set{l \in \bZ: \ \mathsf{x}_l = \mathsf{y}_l},
\]
for some $\gamma > 1$.  Any two such metrics are equivalent.  We can extend such a metric to vertices of $\tree$ by taking the
minimum over all ends that contain the vertices. The set of ends of $\tree$ is the topological boundary of $\tree$ in any of
these metrics.

The dynamics that we consider are maps that preserve the genealogical structure.

\begin{definition}\label{defn: Dyn - abstract}
Let $\tree$ be a tree.  A map $\F: \tree \to \tree$ \emph{preserves children} if for all $\mathsf{v}\in \tree$ the image of a
child of $\mathsf{v}$ is a child of $\F(\mathsf{v}) $.  Symbolically $\F(\vc) = \F(\mathsf{v})^{\child}$.

\end{definition}

A children-preserving map induces a well-defined map on the set of ends of the tree. Additionally, such a map is continuous with
respect to any Gromov metric.  It is easy to check that if $\F: \tree \to \tree$ is a children-preserving map, then there exists
$H \in \bZ$ such that $\F(\tree_l) = \tree_{l-H}$ for all $l \in \bZ$.

\begin{definition}\label{defn: Twd - abstract}
A \emph{tree with dynamics} is a pair $(\tree,\F)$, where $\tree$ is a genealogical tree and $\F: \tree \to \tree$ preserves
children.

\end{definition}

This is a very general definition.  In order to be a tree with dynamics of a polynomial (Definitions \ref{defn: TwD of Poly
w/Disc}, \ref{defn: TwD of Poly w/Con}), there are a number of additional conditions which must be satisfied (see Proposition
\ref{prop: Poly TwD properties} and \cite[Def. 4.7]{E03}).

Throughout this paper, let $(\tree, \F)$ denote a tree with dynamics such that for some $H \in \bZ$, $\F(\tree_l) = \tree_{l-H}$
for all $l$.

We give a few examples of trees with dynamics. Let $(\tree, \F)$ be a tree with dynamics and let $n \geq 1$. A straightforward
check of the above definitions shows that $(\tree, \F^n)$ is also a tree with dynamics.

The following example describes the tree with dynamics of every quadratic polynomial with a disconnected Julia set (see Example
\ref{eg: disc quad TwD}).

\begin{figure}[hbt]
\begin{center}

\includegraphics{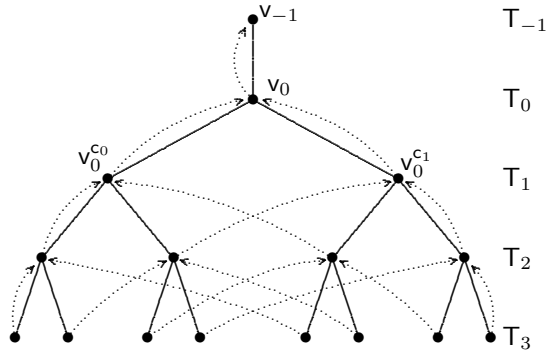}

\caption{A rooted tree with dynamics.  This is the tree with dynamics of every quadratic polynomial with a disconnected Julia set
(see Example \ref{eg: disc quad TwD}).}
\label{fig: Binary tree}
\end{center}
\end{figure}

\begin{example}[Disconnected Quadratic Tree with Dynamics] \label{eg: binary TwD}
We define a tree with dynamics as follows (see Fig.\ \ref{fig: Binary tree}).  Define $\tree_{-l}$ as a single vertex
$\mathsf{v}_{-l}$ for $l \in \bN$. Define $\mathsf{v}_{-l}^{\parent} = \mathsf{v}_{-l-1}$, and $\F(\mathsf{v}_{-l}) =
\mathsf{v}_{-l-1}$  for $l \in \bN$. Give $\mathsf{v}_0$ two children: $\mathsf{v}_0^{\child_0}$ and $\mathsf{v}_0^{\child_1}$,
so $\tree_1 = \set{\mathsf{v}_0^{\child_0}, \mathsf{v}_0^{\child_1}}$. Define $\F(\mathsf{v}_0^{\child_i} ) = \mathsf{v}_0$. For
$\mathsf{v} \in \tree_l $ ($l > 0$), give $\mathsf{v}$ two children $\mathsf{v}^{\child_0}, \mathsf{v}^{\child_1}$, and define
$\F (\mathsf{v}^{\child_i} ) = \F(\mathsf{v})^{\child_i}$. In this tree with dynamics, $H = 1$.
\end{example}

The following example is not natural, but it is a useful source of counter-examples.

\begin{example} \label{eg: Z2 TwD}
Let $\tree = \set{(l,m) \in \bZ^2}$.  Define $(l,m)^{\parent} = (l - 1,\lfloor m/2 \rfloor) $ for $m \geq 0$, and
$(l,m)^{\parent} = (l -1, 0) $ for $m < 0$ (where $\lfloor \cdot \rfloor$ is the greatest integer function). So $\tree_l =
\set{(l,m): \ m\in \bZ}$. For any $H \in \bZ$, define $\F_H(l,m) = (l-H, m)$. Check that $\F_H$ is child preserving (in fact it
is an automorphism of $\tree$). Thus, $(\tree, \F_H)$ is a tree with dynamics. On the other hand, we can define an
infinite-to-one map by $\mathsf{G}_H(l,m) = (l-H, 0)$ for any $H \in \bZ$. Nonetheless $\mathsf{G}_H$ is child preserving, so
$(\tree, \mathsf{G}_H)$ is a tree with dynamics.
\end{example}

The above example is not the tree of any polynomial, since it does not satisfy the first 3 conclusions of Proposition \ref{prop:
Poly TwD properties}.

\subsection{Dynamics of Ends}\label{subs: Twd - Dyn of Ends}

In this paper, the main type of dynamics that we consider are the returns of an end of a tree with dynamics to itself.  Let
$\mathsf{x}_L$ be a vertex of an extended end $\boldsymbol{\mathsf{x}} = (\mathsf{x}_l)$. We say $\F^n(\mathsf{x}_L)$
\emph{returns (to $\boldsymbol{\mathsf{x}}$)} if $\F^n(\mathsf{x}_L) = \mathsf{x}_m$ for some $m \in \bZ$ and $n \geq 1$.  We say
$\F^n(\mathsf{x}_L)$ is the \emph{first return} of $\mathsf{x}_L$ (to $\boldsymbol{\mathsf{x}}$) if $\F^n(\mathsf{x}_L)$ returns
and $n \geq 1 $ is the minimal iterate that returns; we call $n$ the \emph{first return time} of $\mathsf{x}_L$.

\begin{lemma}\label{lem: ancestors return with descendants}
Let $\boldsymbol{\mathsf{x}} = (\mathsf{x}_l)_{l \in \bZ}$ be an extended end of $(\tree, \F)$. Let $\mathsf{x}_L \in
\boldsymbol{\mathsf{x}}$. If $\F^n(\mathsf{x}_L)$ returns for some $n\geq 1$, then $\F^n(\mathsf{x}_{l})$ also returns for every
$l < L$.
\begin{proof}
It suffices to show that $\F^n(\mathsf{x}_{L-1})$ returns. Say that $\F^n(\mathsf{x}_L) = \mathsf{x}_M$. Since $\F$ preserves
children, it also preserves parents.  We have $\F^n(\mathsf{x}_{L-1}) = \F^n(\mathsf{x}_{L}^{\parent}) = \mathsf{x}_{M}^{\parent}
= \mathsf{x}_{M-1}$ by indexing of the end.
\end{proof}
\end{lemma}

\begin{figure}[hbt]

\begin{center}

\includegraphics{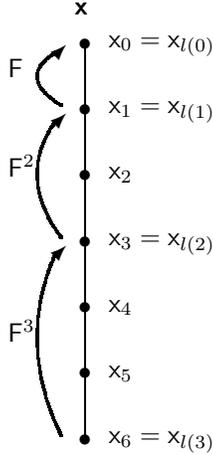}
\caption{An end with a return chain marked.}
 \label{fig: ret chain}
\end{center}

\end{figure}

We introduce a simple combinatorial object, a \emph{return chain} (see Fig.\ \ref{fig: ret chain}). Which is a subset of an end
defined using first returns. It corresponds to a return nest (Definition \ref{defn: return nest}). Essentially this allows us to
consider a ``one-dimensional'' system, rather than the whole tree.

\begin{definition} \label{defs: ret chain}
Let $\boldsymbol{\mathsf{x}}$ be an extended end of $\tree$.  A \emph{return chain} is a sub-end $(\mathsf{x}_{l(k)})_{k \in
\bZ}$ such that
\begin{equation*}
  \F^{n_k} (\mathsf{x}_{l(k)}) = \mathsf{x}_{l(k-1)}
  \quad \text{for all } k \in \bZ,
\end{equation*}
where $n_k$ is the first return time of $\mathsf{x}_{l(k)}$. We call $(n_k)_{k \in \bZ}$ the \emph{return times} of the return
chain. If $l(k) = \min \set{l: \mathsf{x}_{l(k-1)} \text{ is the first return of } \mathsf{x}_l}$ for all $k >0$, then we call
the return chain \emph{minimal}.
\end{definition}

Let $(\tree, \F)$ be a tree with dynamics such that $H \geq 1$, and let $(\mathsf{x}_{l(k)})$ be a return chain.  For each $k$,
$\mathsf{x}_{l(k)}$ is the first return of $\mathsf{x}_{l(k+1)}$. It is possible that some $\mathsf{x}_{l(k)}$ is the first
return of several $\mathsf{x}_{l}$. In a minimal return chain, $\mathsf{x}_{l(k+1)}$ is the first vertex of
$\boldsymbol{\mathsf{x}}$ below $\mathsf{x}_{l(k)}$ whose first return is $\mathsf{x}_{l(k)}$.

\begin{lemma} \label{lem: n_k non-dec}
Let $(\tree, \F)$ be a tree with dynamics such that $H \geq 1$. If $(n_k)_{k \in \bZ}$ are the return times of a return chain
$(\mathsf{x}_{l(k)})_{k \in \bZ}$, then $(n_k)$ is a non-decreasing sequence.
\begin{proof}
Fix $k$. We have $l(k-1) < l(k)$ since $H \geq 1$.  Because $\F^{n_k}(\mathsf{x}_{l(k)})$ is a return,
$\F^{n_k}(\mathsf{x}_{l(k-1)})$ is also a return by Lemma \ref{lem: ancestors return with descendants}.  Thus, $n_k \geq n_{k-1}$
since $n_{k-1}$ is the first return time of $\mathsf{x}_{l(k-1)}$ and first return times are minimal.
\end{proof}

\end{lemma}

Return times are non-decreasing, but we can have $n_{k+1} = n_k$ for some $k$. A \emph{cascade of central returns of length} $j
\geq 1$ is a constant string $n_k= \dots = n_{k+j-1}$ for some $k \geq 1$ (compare to \cite[\S3]{Lyubich-Quadratics}).  The
restriction on $k$ is to avoid trivial cascades of infinite length (see Corollary \ref{cor: n_k = 1 for ret nest}). Cascades of
central returns of infinite length, equivalently eventually constant return times, occur exactly when we have a periodic end.
Thus, we code trivial dynamics by trivial return times.

\begin{lemma}\label{lem: periodic iff n_k = N}

Let $\boldsymbol{\mathsf{x}}$ be an extended end of $(\tree, \F)$.  If $H \geq 1$, then the following are equivalent.
\begin{enumerate}
  \item The minimum period of $\boldsymbol{\mathsf{x}}$ is $N$.

  \item The return times $(n_k)$ of some return chain of $\boldsymbol{\mathsf{x}}$ satisfy $n_k = N$ for all $k$ sufficiently
      large.

  \item   There is a return chain of $\boldsymbol{\mathsf{x}}$, and the return times $(n_k)$ of every return chain of
      $\boldsymbol{\mathsf{x}}$ satisfy $n_k = N$ for all $k$ sufficiently large.

\end{enumerate}

\begin{proof}
$(1 \Rightarrow 2)$.  If $N$ is the minimum period of $\boldsymbol{\mathsf{x}}$, then there is an  $L$ such that for all $l \geq
L$, $\F^n(\mathsf{x}_l) \notin \boldsymbol{\mathsf{x}}$ for $n =1, \dots, N-1$, but $\F^N(\mathsf{x}_l) \in
\boldsymbol{\mathsf{x}}$. So, the first return time of $\mathsf{x}_l$ equals $N$. Define $l(0) = L$.  For $k \leq 0$, define
$n_k$ as the first return time of $\mathsf{x}_{l(k)}$, and $\mathsf{x}_{l(k-1)} = \F^{n_k}(\mathsf{x}_{l(k)})$.  For $k > 0$,
define $n_k = N$, and $l(k) = l(k-1)+HN$.  Then $(\mathsf{x}_{l(k)})_{k \in \bZ}$ is a return chain with $n_k = N$ for all $k
\geq 0$.

$(2 \Rightarrow 3)$. Let $(\mathsf{x}_{l(k)})_{k \in \bZ}$ be a return chain with $n_k = N$ for all $k \geq K$.  Let
$(\mathsf{x}_{\hat{l}(k)})_{k \in \bZ}$ be a return chain of $\boldsymbol{\mathsf{x}}$ with return times $\hat{n}_k$.   Fix
$\hat{l}(k) \geq l(K)$. Then $\hat{n}_k \geq N$ by Lemma \ref{lem: ancestors return with descendants}. Also, we can find $K' > K$
such that $ l(K') > \hat{l}(k)$. So, $\hat{n}_k = N$ by Lemma \ref{lem: ancestors return with descendants}.

$(3 \Rightarrow 1)$. Let $(\mathsf{x}_{l(k)})_{k \in \bZ}$ be a return chain of $\boldsymbol{\mathsf{x}}$.  Then the return times
satisfy $n_k = N$ for all $k \geq K$ for some $K$. For any $l \in \bZ$, we can find $k \geq K$ such that $l(k) > l$. Since $n_k =
N$, we have $\F^N(x_l) \in \boldsymbol{\mathsf{x}}$ by Lemma \ref{lem: ancestors return with descendants}. Therefore
$\F^N(\boldsymbol{\mathsf{x}}) = \boldsymbol{\mathsf{x}}$, and $\boldsymbol{\mathsf{x}}$ is periodic with period $N$.  For $l
\geq l(K)$, $F^n(\mathsf{x}_l) \notin \boldsymbol{\mathsf{x}}$ for $n = 1, \dots, N-1$ by Lemma \ref{lem: ancestors return with
descendants} and minimality of $n_K$.  Therefore, $N$ is the minimum period of $\boldsymbol{\mathsf{x}}$.

\end{proof}

\end{lemma}

\subsection{Main Theorem for Trees with Dynamics}\label{subs: TwD - Main}

We are ready prove our main theorem for trees with dynamics, which is a combinatorial version of Theorem \ref{thm: Main}.

\begin{theorem}\label{thm: ret times of ret chains Fib}
Let $(\tree, \F)$ be a tree with dynamics such that $H \geq 1$.  The return times of any return chain of $(\tree, \F)$ are
variable-$r$ meta-Fibonacci numbers.
\begin{proof}
Let $(\mathsf{x}_{l(k)})_{k \in \bZ}$ be a return chain with return times $(n_k)_{k \in \bZ}$. Fix $k \in \bZ$. If $n_{k} =
n_{k-1}$, then $r(k) = 1$ and we are done. Otherwise, $n_k > n_{k-1}$ by Lemma \ref{lem: n_k non-dec}. Thus, $n_{k} = N +
n_{k-1}$ for some $N \in \bZ^+$. We have
\[
    \F^{n_k}(\mathsf{x}_{l(k-1)}) = \F^{N + n_{k-1}}(\mathsf{x}_{l(k-1)})
    = \F^{N }(\F^{n_{k-1}}(\mathsf{x}_{l(k-1)})) = \F^N(\mathsf{x}_{l(k-2)}).
\]
Since $\F^{n_k}(\mathsf{x}_{l(k)})$  returns, so does $\F^{n_k}(\mathsf{x}_{l(k-1)})$ by Lemma \ref{lem: ancestors return with
descendants}.  Thus, $N \geq n_{k-2}$ by minimality of $n_{k-2}$.  Therefore, $n_k \geq n_{k-1} + n_{k-2}$.  If we have equality,
then we are done.  Otherwise, we can repeat the above argument to show that $n_k \geq n_{k-1} + n_{k-2} + n_{k-3}$ and so on.
After a finite number of repetitions of this argument we will have equality, since $n_k$ is finite. Therefore for some $r(k) \geq
1$,
\[
    n_k = n_{k-1} + \dots + n_{k-r(k)}.
\]
\end{proof}
\end{theorem}

There is no reason that we should expect that $r(k) \leq k$. Hence, we must consider \emph{extended} $r(k)$-bonacci numbers.  In
fact given $K, R \in \bZ^+$, there is a tree with dynamics of a polynomial with a disconnected Julia set, which has a return
chain with a return times $(n_k)$ generated by $r(k)$ such that $r(K) \geq R$ \cite[Lem.\ 7.12]{E03}.

Our main theorem follows from the above theorem.  An end of a tree with dynamics of a polynomial corresponds to a sequence of
nested puzzle pieces (Definition \ref{defn: nest of puzzle pieces}).  Thus, a return chain corresponds to a return nest
(Definition \ref{defn: return nest}). We need only compare Definitions  \ref{defs: ret chain} and \ref{defn: return nest}, and
verify that they define the same return times.  The details of the correspondence are explained in Section \ref{sect: poly TwD}.

\subsection{Constructing Return Chains}\label{subs: TwD - Ret Chain}

Say an end is \emph{recurrent} if its forward iterates accumulate at itself. First, we show that every vertex of a recurrent end
is the first return of some vertex of the end. Notice by definition of the metric on trees, an end $\boldsymbol{\mathsf{x}}$ is
recurrent if and only if for all $L \in \bZ$ there is $K \in \bZ$ such that $\F^N(\mathsf{x}_K) = \mathsf{x}_L$ for some $N \geq
1$.

\begin{lemma} \label{lem: 1 hit}
Let $(\tree, \F)$ be a tree with dynamics.  Let $\boldsymbol{\mathsf{x}}$ be an extended end of $\tree$ that is recurrent under
$\F$. For every $L \in \bZ$, there exists $M \in \bZ$ such that $\mathsf{x}_L$ is the first return of $\mathsf{x}_M$.

\begin{proof}
Fix $L \in \bZ$. Since $\boldsymbol{\mathsf{x}}$ is recurrent, $\F^n(\mathsf{x}_K) = \mathsf{x}_L$ for some $K \in \bZ$ and $n
\geq 1$.  It is possible that $\mathsf{x}_K$ has returned to $\boldsymbol{\mathsf{x}}$ before it returns to $\mathsf{x}_L$. Let
$F^{N}(\mathsf{x}_K) = \mathsf{x}_M$ be the last return of $\mathsf{x}_K$ before it returns to  $\mathsf{x}_L$. Then $F^{n -
N}(\mathsf{x}_M) = \mathsf{x}_L$ and this is the first return of $\mathsf{x}_M$.

 \end{proof}

\end{lemma}

We can form a return chain of an end if and only if the end is recurrent.

\begin{prop} \label{prop: form min ret chain}
Let $(\tree, \F)$ be a tree with dynamics with $H \geq 1$.  Let $\boldsymbol{\mathsf{x}}$ be an extended end of $\tree$ that is
recurrent under $\F$. For any $L \in \bZ$,  $\boldsymbol{\mathsf{x}}$ has a unique minimal return chain
$(\mathsf{x}_{l(k)})_{k\in \bZ}$ with $l(0)= L$.

\begin{proof}
Let $l(0) =L$.  For $k \leq 0$, recursively define $n_k$ as the first return time of $\mathsf{x}_{l(k)}$, and $l(k-1)$ by
$\mathsf{x}_{l(k-1)} = \F^{n_k}(\mathsf{x}_{l(k)})$. Notice there is no reason to believe that $l(k)$ is minimal for $k \leq 0$,
nor does Definition \ref{defs: ret chain} required it.  For $k >0$,  $\mathsf{x}_{l(k)}$ is the first return of some vertex of
$\boldsymbol{\mathsf{x}}$ by Lemma \ref{lem: 1 hit}.   There is such a vertex with least index since $H \geq 1$.  We define
$l(k+1)$ as the least integer such that $\mathsf{x}_{l(k)}$ is the first return of $\mathsf{x}_{l(k+1)}$.  Define $n_{k+1}$ as
the first return time of $\mathsf{x}_{l(k+1)}$.  We have made no choices, so the minimal return chain we have constructed is
unique.
\end{proof}

\end{prop}

If $\boldsymbol{\mathsf{x}}$ is not recurrent, then the above construction will break down at some stage. That is for some  $K
\geq 0$, $\mathsf{x}_{l(K)}$ will not be the first return of any vertex of $\boldsymbol{\mathsf{x}}$.  In this case, we define
$n_k = r(k) = \infty$ for all $k > K$.  The above construction is a combinatorial version of the construction of a principal nest
\cite[$\S$3.1]{Lyubich-Quadratics}.

Let $(\mathsf{x}_{l(k)})$ be a return chain of $\boldsymbol{\mathsf{x}}$ with return times $(n_k)$ generated by $r(k)$.  By Lemma
\ref{lem: periodic iff n_k = N}, $\boldsymbol{\mathsf{x}}$ is periodic if and only if $r(k) = 1$ for all $k$ sufficiently large.
From the above construction, $\boldsymbol{\mathsf{x}}$ is non-recurrent if and only if $r(k) = \infty$ for all $k$ sufficiently
large.  Therefore, $r(k)$ can be thought of as an indication of the degree to which $\boldsymbol{\mathsf{x}}$ is recurrent: Small
$r(k)$ means that $\boldsymbol{\mathsf{x}}$ is highly recurrent; large $r(k)$ means that $\boldsymbol{\mathsf{x}}$ is weakly
recurrent.

\subsection{Rooted Trees}\label{subs: Twd - Rooted Trees}

We have been working with quite general trees with dynamics.  We now restrict our attention to trees with a distinguished vertex,
a \emph{root}. The tree with dynamics of any polynomial has a root (Proposition \ref{prop: Poly TwD properties}.3).

\begin{definition}
Let $ \tree$ be a tree, and let $\mathsf{v}_0 \in \tree$.  We call $\mathsf{v}_0$ the \emph{root} of $\tree$, if every ancestor
of $\mathsf{v}_0$ has exactly one child and $\mathsf{v}_0$ has more than one child. We call a tree with a root a \emph{rooted
tree}.
\end{definition}

The root of a tree is unique, if it exists.  In a rooted tree, we index the levels so that $\tree_0 = \set{\mathsf{v}_0}$, where
$\mathsf{v}_0$ is the root of $\tree$. The ancestors of the root are a countable line graph. Thus, $\tree_{-l} $ is a single
vertex $\mathsf{v}_{-l}$ for $l \in \bN$.   The boundary of a rooted tree has the topology of a Cantor set (compact, perfect, and
totally disconnected) union one point (corresponding to $\lim \mathsf{v}_{-l}$).  The dynamics of a rooted tree are Lipschitz
continuous. There are restrictions on return times in a rooted tree.

\begin{prop}\label{prop: n_k = 1 for k < K}
Let $(\tree, \F)$ be a tree with dynamics such that $H \geq 1$. Let $(\mathsf{x}_{l(k)})_{k \in \bZ}$ be a return chain of
$\boldsymbol{\mathsf{x}}$ with return times $(n_k)$ generated by $r(k)$. If $\tree$ is rooted, then there exists $K \in \bZ$ such
that $n_k = 1$ and $r(k) = 1$ for all $k \leq K$.

\begin{proof}
Since $\tree$ is rooted, we have $\tree_{-l} = \set{\mathsf{v}_{-l}}$ for all $l \in \bN$.  Thus, $\mathsf{x}_{-l} =
\mathsf{v}_{-l}$, and $\F(\mathsf{x}_{-l}) = \mathsf{x}_{-l -H} = \mathsf{v}_{-l -H} $ for all $l \in \bN$. We can find $K$ such
that $l(K) \leq H$. Then $\F(\mathsf{x}_{l(K)}) = \mathsf{x}_{l(K)-H}$, so $n_K =1$. Since return times are non-decreasing by
Lemma \ref{lem: n_k non-dec}, $n_k = 1$ for all $k \leq K$. By definition, $r(k) = 1$ for all $k \leq K$.
\end{proof}
\end{prop}

This proposition  describes the initial conditions satisfied by the return times of a return chain in a rooted tree with
dynamics. When we construct a return chain, we can choose $l(0)$ (Proposition \ref{prop: form min ret chain}). We define a normal
form for return chains by requiring that $l(0) = 0$.  This requirement normalizes the initial conditions of the return times.

\begin{cor}\label{cor: n_k = 1 PL chain}
Let $(\tree, \F)$ be a rooted tree with dynamics such that $H \geq 1$.  Let $(\mathsf{x}_{l(k)})_{k \in \bZ}$ be a  return chain
of $\boldsymbol{\mathsf{x}}$ with return times $(n_k)$ generated by $r(k)$.  If  $l(0) = 0$, then $n_k = r(k) =1$ for all $k \leq
0$.  Furthermore, if $(\mathsf{x}_{l(k)})$ is minimal and $l(0) = 0$, then $n_1 = r(1) =1$.

\begin{proof}
Since $l(0) = 0 < H$, it follows from the proof of Proposition \ref{prop: n_k = 1 for k < K} that $n_k = 1$ for all $k \leq 0$.

Now, assume that $(\mathsf{x}_{l(k)})$ is minimal. We have $\F(\mathsf{x}_{H}) = \mathsf{x}_{0} = \mathsf{x}_{l(0)}$.  Thus $l(1)
= H$, since $H$ is the least possible index of a vertex of $\boldsymbol{\mathsf{x}}$ that could return to $\mathsf{x}_{0} $.
Therefore, $n_1 = r(1) =1$.
\end{proof}
\end{cor}

\section{Variable-$r$ Meta-Fibonacci Sequences} \label{sect: meta-Fib seqs}

In this section, we study variable-$r$ meta-Fibonacci sequences. In light of Lemma \ref{lem: n_k non-dec} and Corollary \ref{cor:
n_k = 1 PL chain}, we will only consider non-decreasing $r(k)$-bonacci sequences $(n_k )$ with $n_k = 1$ for $k \leq 0$. We
derive some estimates for $n_k$ based on bounds of $r(k)$ (\S\ref{subs: MF - Est}). We give two main estimates for $n_k$, a lower
bound and  an upper bound (Propositions \ref{prop: Lower Bound} and \ref{prop: Upper Bound} respectively). We also give some
examples of variable-$r$ meta-Fibonacci sequences (\S\ref{subs: MF - Egs}).

We say an $r(k)$-bonacci sequence $(n_k)$ has a \emph{cascade} of length $j \geq 1$ if $n_k =  \dots=  n_{k+j-1}$ for some $k
\geq 1 $. Naturally, if $(n_k)$ is the sequence of return times of a return chain, then a cascade corresponds to a cascade of
central returns.

Let $\lfloor \cdot \rfloor$ denote the greatest integer function.

\begin{prop}\label{prop: Lower Bound}
Let $(n_k)_{k \in \bZ}$ be a non-decreasing $r(k)$-bonacci sequence such that $n_k = 1$ for all $k \leq 0$. Let $J \in \bZ^+$. If
the length of every cascade of $(n_k)$ is bounded above by $J$, then
\[
    n_k \geq 2^{\lfloor k/(J+2) \rfloor} \quad \text{for every $k \geq 1$}.
\]

\end{prop}

\begin{prop}\label{prop: Upper Bound}
Let $(n_k)_{k \in \bZ}$ be a non-decreasing $r(k)$-bonacci sequence such that $n_k = 1$ for all $k \leq 0$.  Let $M \in \bZ^+$.
If $r(k+1) \leq M r(k) + 1$ for every $k\geq 1$, then
\[
    n_k \leq (M+1)^k \quad \text{for every $k \geq 1$}.
\]
\end{prop}

\subsection{Examples}\label{subs: MF - Egs}

Before we prove the above estimates, we give some examples of $r(k)$-bonacci sequences. The most elementary example is when
$r(1)= 1$ and $r(k) = 2$ for all $k \geq 2$. Then $n_k = u_{k+1}$, where $(u_k)$ is the usual Fibonacci sequence. Suppose that
$r(k) = r \geq 2$ for all $k$ sufficiently large.  Then up to re-indexing, the tail of $n_k$ is a generalized $r$-bonacci
sequence for some initial conditions \cite{Miles}.

If $r(k) =1$ for all $k$, then $n_k=1 $ for all $k$.  If $r(k) =1$ for all $k$ large, then $n_k$ is eventually constant.  While
this is may seem like a trivial example, it corresponds to the return times of a periodic end (Lemma \ref{lem: periodic iff n_k =
N}).

The following example shows that the return times of the Feigenbaum polynomial \cite{Sullivan_renormalization} are $r(k)$-bonacci
numbers.

\begin{example} \label{eg: n_k = 2^k-1}
Let $r(k) = k$ for all $k \geq 1$.  By an easy inductive argument, we find that $n_k = 2^{k-1}$ for $k \geq 1$.
\end{example}

We can make $n_k$ grow linearly by taking $r(k) =1 $ for many successive $k$.

\begin{example} \label{eg: lin growth}
For $k \geq 2$, let $r(k) = 2$ if $k = 2^m$ for some $m \in \bZ^+$, and let $r(k) = 1$ otherwise.
\[
\begin{array}{|c|c|c|c|c|c|c|c|c|c|c|c|c|}\hline
  k    & 0 & 1  &  2 &  3 & 4  & 5  &  6 & 7  & 8  & 9  & \cdots & 16\\
  \hline
  r(k) & 1 & 1  &  2 &  1 & 2  & 1  &  1 & 1  & 2  & 1  & \cdots & 2\\
  \hline
  n_k  & 1 & 1  &  2 &  2 &  4 & 4  & 4  & 4  & 8  & 8  & \cdots & 16 \\
  \hline
\end{array}
\]
It is straightforward to show that $k/2 < n_k \leq k$ for $k \geq 1$.
\end{example}

The following is an example of an \emph{extended} $r(k)$-bonacci sequence.

\begin{example}  \label{eg: r(k) = 2^k}
Define $r(k) = 2^{k-1}$ for $k \geq 1$.
\[
\begin{array}{|c|c|c|c|c|c|c|c|c|c|c|c|}
\hline
  k    & \leq 0 & 1  &  2 &  3 & 4   & 5    &  6   & 7    & 8    & 9 & 10 \\ \hline
  r(k) & 1      & 1    & 2  &  4 & 8  &  16 &  32  & 64   & 128  & 256 &  512 \\
  \hline
  n_k  & 1      & 1 & 2  &  5 & 13  &  33 &  81 & 193  & 449  & 1025  & 2305  \\
  \hline
\end{array}
\]
It follows from Lemma \ref{lem: Delta r(k) = 1 -> growth =2} and Proposition \ref{prop: Upper Bound} that $2 < n_k/n_{k-1} < 3$
for every $k \geq 3$.
\end{example}

\subsection{Estimates}\label{subs: MF - Est}

We give some estimates on the terms of an $r(k)$-bonacci sequence based on bounds on $r(k)$. For a given $k$, the larger $r(k)$
is, the larger $n_k$ will be in comparison to $n_{k-1}$.  However the exact relationship is subtle.  We are particularly
interested in closed-form bounds. That is, bounds for $n_k$ which are functions of $k$, but not of $n_1, \dots, n_{k-1}$ or
$r(k)$.

First, we recall the asymptotic growth rate of the $r$-bonacci numbers. Let $r \in \bZ^+ $. Define $\gamma_r$ as the unique root
of
\[
    z^r - z^{r-1} - \dots - z - 1
\]
such that $1 \leq |\gamma_r| < 2$, where $|\cdot|$ denotes the complex norm.  It is known that $\gamma_r$ is well-defined and
real for all $r$ \cite[Eq.\ 6$''$]{Miles}. For example, $\gamma_2 = (1+\sqrt{5})/2$.  Let $(u_{r,k})_{k=1}^{\infty}$ be the
$r$-bonacci numbers, then $ \lim_{k \to \infty } u_{r, k+1}/u_{r, k} = \gamma_r$. The sequence $(\gamma_r)_{r=1}^{\infty}$ is
strictly increasing, and $\lim_{r \to \infty} \gamma_r$ = 2 \cite{Dubeau}. We can compare the growth rate of an $r(k)$-bonacci
sequence with bounded $r(k) $ to $\gamma_r$.

Let $(n_k)_{k \in \bZ}$ be a non-decreasing $r(k)$-bonacci sequence such that $n_k = 1$ for all $k \leq 0$. Let $  R \in \bZ^+$.
\begin{equation}
    \text{If } \liminf_{k \to \infty } r(k) \geq R, \quad \text{then }
   n_k \geq \textnormal{const.}\, \gamma_R^k \quad \text{for all $k \geq 1$}
\end{equation}
and some positive constant \cite[Prop.\ 3.5]{E06}.  A similar upper bound holds.
\begin{equation}
    \text{If } \limsup_{k \to \infty } r(k) \leq R, \quad \text{then }
   n_k \leq \textnormal{const.}\, \gamma_R^k \quad \text{for all $k \geq 1$}
\end{equation}
and some positive constant.

We now derive upper and lower bounds using weaker assumptions.  The following lemma is the key observation about the growth of
$n_k$ as a function of $r(k)$. We give a condition for $n_k$ to double (see Example \ref{eg: n_k = 2^k-1}).

\begin{lemma} \cite[Lem.\ 3.3]{E06} \label{lem: Delta r(k) = 1 -> growth =2}
Let $(n_k)_{k \in \bZ}$ be an $r(k)$-bonacci sequence.   If $r(k+1) = r(k) + 1$ for some $k$, then $n_{k+1} =2n_k$.

\begin{proof}
By definition,
\[
     n_{k+1} = \sum_{j=1}^{r(k+1)} n_{k+1-j}
        = n_{k} + \sum_{j=2}^{r(k)+1} n_{k + 1 -j}
        =  n_{k} + \sum_{i=1}^{r(k)} n_{k-i}
        = 2n_{k}.
\]
\end{proof}
\end{lemma}

There are $r(k)$-bonacci sequences where $(n_k)$ grows linearly (Example \ref{eg: lin growth}), or even slower \cite[Thm.\
1]{E06}.  This slow growth occurs when we have long cascades. Proposition \ref{prop: Lower Bound} shows that this is the only way
to obtain sub-exponential growth of $n_k$.

\begin{proof}[Proof of Proposition \ref{prop: Lower Bound}]
First, consider the assumption that the length of cascades of $(n_k)$ is bounded above by $J$. It follows that for $k \geq 1$,
the maximum number of consecutive $r(k)$ that equal 1 is also $J$. Hence, at least one of  $r(k+1), \dots, r(k+J+1)$ is larger
than 1 for any $k \geq 0$.

Next fix $k \geq 0$. It suffices to show that $n_{k+J+2} \geq 2 n_k$.  If $r(k+1)= 1$ or $r(k+2) = 1$, then $n_{k+J+2} \geq 2
n_k$ by the above observation about $r(k)$ and Lemma \ref{lem: Delta r(k) = 1 -> growth =2}.  By similar arguments, we can reduce
to the case when $r(k+1) \geq r(k+2) > 1$.  This implies that $n_{k+1} \geq n_k  + n_{k-1}$, and $n_{k+2} \geq n_{k+1}  +n_k $.
Therefore, $ n_{k+J+2} \geq n_{k+2} \geq 2n_{k} + n_{k-1} > 2n_{k}$.

\end{proof}

The above estimate is sharp, as the following example shows. Fix $J \geq 1$. Define $r(k) = 2$ if  $k \equiv 0 \ (\Mod J+2)$ and
$k\geq 1$, and $r(k) =1$ otherwise. Then for $k\geq 1$, $n_k = 2^{k/(J+2)}$ if $k \equiv 0 \ (\Mod J+2)$ and $k\geq 1$, and $n_k
= n_{k-1}$ otherwise.

Proposition \ref{prop: Upper Bound} follows from the following lemma by induction.  We derive an upper bound for the magnitude of
$n_{k+1}$ relative to $n_k$.

\begin{lemma}
Let $(n_k)_{k \in \bZ}$ be a non-decreasing $r(k)$-bonacci sequence with $n_k =1$ for $k\leq 0$. Let $M \in \bZ^+ $. If $r(k+1)
\leq Mr(k) + 1 $ for some $k$, then $n_{k+1} \leq (M+ 1)n_{k}$. Moreover, the equivalent statement with strict inequalities
holds.

\begin{proof}
By definition,
\begin{align*}
     n_{k+1} &= n_{k} + \sum_{j=2}^{r(k+1)} n_{k + 1 -j}\\
        &\leq  n_{k} + \sum_{i=1}^{Mr(k)} n_{k-i} \quad \text{(by assumption)} \\
        &=  n_{k} +  \sum_{m=0}^{M-1} \sum_{i=1}^{r(k)} n_{k-i- mr(k)}\\
        &\leq n_{k} +  \sum_{m=0}^{M-1} \sum_{i=1}^{r(k)} n_{k-i} \quad \text{(since the $n_k$ are non-increasing)}\\
        &= (1+M)n_k.
\end{align*}
Moreover, if $r(k+1) < Mr(k) + 1 $, then the second line above is a strict inequality.
\end{proof}

\end{lemma}

The above estimate is sharp.  If $n_{k-M}= \dots= n_{k}$ and $r(k+1) =M+1$ for some $k$, then $r(k+1) = Mr(k) + 1$ and $n_{k+1} =
(M+1)n_k$.

The above upper bound is unexpectedly small.  A constant function $r(k)$ generates a sequence $(n_k)$ which grows exponentially.
A logarithmic $r(k)$ generates a sequence which grows linearly (Example \ref{eg: lin growth}).  A priori, if $r(k)$ grows
exponentially, then we might expect that $(n_k)$ would grow super-exponentially. However, Proposition \ref{prop: Upper Bound}
shows that it does not (see Example \ref{eg: r(k) = 2^k}).

\section{ Polynomial Trees with Dynamics} \label{sect: poly TwD}

In this section, we describe the construction of a tree with dynamics associated with a polynomial.  We also prove our main
theorem. Every polynomial of degree at least two is associated with a tree with dynamics.  First, we define abstract puzzle of a
polynomial, which is a sequence of decompositions of the Julia set of the polynomial (\S\ref{subs: Poly - Puzzles}).  Every
puzzle has a tree structure. A function that respects the puzzle structure gives rise to a tree with dynamics.  We prove Theorem
\ref{thm: Main}. A polynomial of degree at least 2 has a standard puzzle. Green's function of a polynomial decomposes the plane
into a puzzle, and the dynamics of the polynomial are compatible with the puzzle structure. We outline the construction of the
Branner-Hubbard puzzle for a polynomial with a disconnected Julia set (\S\ref{subs: Poly - Disc J}). This puzzle defines a
canonical tree with dynamics of a polynomial with a disconnected Julia set.  A polynomial with a connected Julia set has a tree
with dynamics for each of its Yoccoz puzzles (\S\ref{subs: Poly - Con J}). We note some common properties of standard polynomial
tree with dynamics (\S\ref{subs: Poly - Poly TwD}).  Finally, we describe a generalized principal nest.

\subsection{Puzzles}\label{subs: Poly - Puzzles}

A puzzle is a sequence of Markov partitions in a general sense.  We define an abstract puzzle of a polynomial. We have two main
reasons for doing so.  First, we wish to list common properties of polynomial puzzles. Second, we wish to isolate the properties
we need to define a tree with dynamics.

\begin{definition}\label{defn: puzzle}
Let $f$ be a polynomial of degree at least 2.   A \emph{puzzle} of $f$ is a sequence $\puz = (\puz_l)_{l \in \bZ}$, where each
$\puz_l$ is an at most countable collection of disjoint non-empty subsets of $\bC$ such that for each $l \in \bZ$,
\[
    \Kjulia(f) \subset \bigcup_{P \in \puz_l} \overline{P}.
\]
We call each $P \in \puz_l$ a \emph{puzzle piece} at depth $l$. We require that the puzzle pieces satisfy the following
\emph{Markov properties}:

\begin{enumerate}
  \item For any puzzle pieces $P_1, P_2$, there is a puzzle piece $P$ such that $P_1 \cup P_2 \subset P$.
  \item If $P \in \puz_l$ for some $l$, then $P \subsetneq P^{\parent}$ for some (unique) $ P^{\parent} \in \puz_{l-1}$.
  \item If $P \in \puz_l$ for some $l$, then $f(P) = Q$ and $f(P^{\parent}) = Q^{\parent}$, for some $Q \in \puz_{L}$ with $L
      < l$.
\end{enumerate}
\end{definition}

There is a standard puzzle for each polynomial (of degree $\geq 2$), which is the Branner-Hubbard puzzle in the disconnected case
\cite{BH92, E03}, or a Yoccoz puzzle in the connected case \cite{Hubbard_Loc_Con, Kiwi_Rat_Rays}.  However, the standard
definition of a Yoccoz puzzle only include $\puz_l$ with $l \in \bN$.  So we modify the standard definition to define $\puz_l$
for $l < 0$ (Definition \ref{defn: TwD of Poly w/Con}).  These modifications produce a finite number of puzzle pieces with
\emph{exceptional dynamics}; they only satisfy a weakened version of condition 3:

\begin{equation*}
\emph{(3') If $P \in \puz_l$ for some $l$, then $f(P) \subset Q \subset f(P^{\parent})$ for some $Q \in \puz_{L} $ with $L < l$.}
\end{equation*}

\begin{definition}\label{dfen: T of puz}
Let $\puz = (\puz_l)_{l \in \bZ}$ be a puzzle of a polynomial $f$. We define $\tree(\puz)$, the genealogical \emph{tree of}
$\puz$, by defining each puzzle piece as a vertex. We define the \emph{parent} of $P$ as $P^{\parent}$.

\end{definition}

We define the dynamics on a tree of a puzzle.   Notice that the above definition defines the \emph{ancestors} of a puzzle piece.

\begin{definition} \label{defn: induced dyn on puz}
Let $\puz = (\puz_l)_{l \in \bZ}$ be a puzzle a polynomial $f$. We define the \emph{induced dynamics} $\F_f: \tree(\puz) \to
\tree(\puz)$ by $\F_f(P) = f(P)$ for non-exceptional puzzle pieces.

Assume that every puzzle piece with exceptional dynamics has an ancestor with non-exceptional dynamics. For a puzzle piece $P$
with exceptional dynamics, inductively assume that $\F_f(P^{\parent})$ is defined. Define $\F_f(P) = Q$, where $Q$ is the unique
puzzle piece such that $f(P) \subset Q$ and $Q^{\parent} = \F_f(P^{\parent})$.

\end{definition}

We only used the Markov properties of the puzzle for the above definitions: puzzle pieces are ordered by inclusion, and $f$
respects this order.

\begin{lemma}\label{lem: TwD of Puz}
Let $\puz$ be a puzzle of a polynomial $f$.  If every puzzle piece with exceptional dynamics has an ancestor with non-exceptional
dynamics, then $(\tree(\puz), \F_f)$ is a tree with dynamics.

\begin{proof}
Definitions \ref{dfen: T of puz}  and \ref{defn: induced dyn on puz} certainly define a graph with a self map.  Markov property 1
implies that this graph is connected, and property 2 implies it has no non-trivial circuits.  Hence, it is a tree.

Property 3 implies that $\F_f$ preserves children.  We must show that $\F_f$ is well defined on puzzle pieces with exceptional
dynamics. Suppose $P$ has exceptional dynamics. Inductively we may assume that we have defined $\F_f(P^{\parent})$, such that
$f(P^{\parent}) \subset \F_f(P^{\parent}) \in\puz_{L+1}$ for some $L$. It follows from Condition 3$'$ that $f(P) \subset Q_2
\subset f(P^{\parent})$ for some puzzle piece $Q_2 \in \puz_{l_2} $ with $l_2 \leq L$. We can find $Q_1 \in \puz_{L}$ such that
$f(P) \subset Q_2 \subset Q _1$ and $Q_1^{\parent} = \F_f(P^\parent)$, by applying Condition 2 a finite number of times. If
$Q_1'$ is another such piece, then $f(P) \subset Q_1 \cap Q_1'$. So $Q_1 = Q_1'$ by disjointness of puzzle pieces at level $L$.
Thus $Q_1 $ is unique, and $\F_f(P) = Q_1$ is well-defined.
\end{proof}
\end{lemma}

 Given a puzzle, we define a nest of puzzle pieces, which corresponds to an end of the tree with dynamics.

\begin{definition}\label{defn: nest of puzzle pieces}
Let $\puz$ be a puzzle of a polynomial $f$.  A \emph{nest of puzzle pieces of $f$} is a sequence $(P_l)_{l \in \bN}$ such that
each $P_l \in \puz_l$ and $P_{l} = P_{l+1}^{\parent}$. An \emph{extended nest} is the analogous sequence with $l \in \bZ$.
\end{definition}

Let $(P_l)$ be an extended nest of puzzle pieces of some puzzle $\puz$ of a polynomial $f$. Let $(\mathsf{x}_l)$ be the extended
end in the tree with dynamics $(\tree, \F) = (\tree(\puz), \F_f)$ that corresponds to $(P_l)$. Note that it is possible that
$\F^n(\mathsf{x}_l) = \mathsf{x}_m$, but  $f^n(P_l) \subsetneq P_m$  for some indices. For $l \in \bZ$ and $n \geq 1$, we say
$f^n(P_l)$ \emph{returns} if $\F^n(\mathsf{x}_l)$ returns.  We define the first return time of $P_l$ as the first return time of
$\mathsf{x}_l$.

\begin{definition}\label{defn: return nest}
Let $(P_l)_{l \in \bZ}$ be an extended nest of puzzle pieces of a polynomial $f$. Let $(\tree(\puz), \F_f)$ be the tree with
dynamics of $\puz$. A \emph{return nest} is a sub-nest $(P_{l(k)})_{k \in \bZ}$ such that
\[
    \F_f^{n_k}(P_{l(k)}) = P_{l(k-1)},
\]
where $n_k $ is the first return time of $P_{l(k)}$.  We call $(n_k)_{k \in \bZ}$ the \emph{return times} of the return nest. If
$l(k) = \min \set{l: P_{l(k-1)} \text{ is the first return of } P_l}$ for all $k > 0$, then we call the return nest
\emph{minimal}.
\end{definition}

If  $f^n(P_{l(k)})$ has exceptional dynamics for for some $k$ and some $n=0, \dots, n_k - 1$, then we only have that
$f^{n_k}(P_{l(k)}) \subset P_{l(k-1)}$. Otherwise, we have $f^{n_k}(P_{l(k)}) = P_{l(k-1)}$. A nest of standard puzzle of $f$ has
at most one exceptional puzzle piece. In a Yoccoz puzzle, if $l(k-1) =0 $, then we may have $f^{n_k}(P_{l(k)}) \subsetneq
P_{l(k-1)}$. A principal nest (Definition \ref{defn: principal nest}) is a special case of a return nest.

\begin{lemma} \label{lem: dictionary}

Let $f$ be a polynomial of degree $\geq 2$. Let $\puz$ be a puzzle of $f$ and $(\tree, \F) = (\tree(\puz), \F_f)$. The following
chart is a dictionary between the dynamical systems $(\puz, f)$ and $(\tree, \F)$:

\begin{figure}[hbt]
\begin{center}
\begin{tabular}{|l|l|}
  \hline
  $(\puz, f)$ & $(\tree, \F)$ \\
  \hline
  $P$ a puzzle piece & $\mathsf{v}$ a vertex \\
  \hline
  $(P_l)_{l \in \bZ}$ an extended nest of puzzle pieces
  & $(\mathsf{x}_l)_{l \in \bZ}$ an extended end \\
  \hline
  $(P_{l(k)})_{k \in \bZ}$ a return nest
  & $(\mathsf{x}_{l(k)})_{k \in \bZ}$ a return chain \\
  with return times $(n_k)$ & with return times $(n_k)$\\
  \hline
\end{tabular}

\end{center}
\end{figure}

\begin{proof}
The only non-obvious statement is that the return times of a return nest and the corresponding return chain are equal. Observe
that the definition of return nest (Definition \ref{defn: return nest}) uses the tree with dynamics.  Therefore, the return times
are the same.

\end{proof}
\end{lemma}

We prove our main theorem: The return times of any return nest are variable-$r$ meta-Fibonacci numbers.

\begin{proof}[\textbf{Proof of Theorem \ref{thm: Main}}]
Theorem \ref{thm: ret times of ret chains Fib} applies to any tree with dynamics of a polynomial by Proposition \ref{prop: Poly
TwD properties}. By the dictionary (Lemma \ref{lem: dictionary}), a return nest corresponds to a return chain and the return
times are the same.
\end{proof}

When we construct a standard puzzle of a polynomial $\puz$, then $\puz_l$ always contains exactly one puzzle piece for $l \leq 0$
(Proposition \ref{prop: Poly TwD properties}). Hence we can say more about the return times of a nest from a standard puzzle.

\begin{cor}\label{cor: n_k = 1 for ret nest}
Let $\puz$ be a puzzle of $f$ such that $\puz_l$ contains exactly on puzzle piece for $l\leq 0$.  If $(n_k)$ are the return times
of a return nest $(P_{l{(k)}})$, then $n_k = r(k) = 1$ for all $k \leq K$ for some $K\in \bZ$. If  $l(0) = 0$, then $K =0$.
Furthermore, if $l(0) = 0$ and the nest is minimal, then $K =1$.
\begin{proof}
Say $\puz_0 = \set{P_0}$.  Since $\puz_l$ contains exactly on puzzle piece for $l\leq 0$, $P_0$ is the root of $(\tree({\puz}),
\F_f)$.  Apply Proposition \ref{prop: n_k = 1 for k < K}.  For the last parts use Corollary \ref{cor: n_k = 1 PL chain}.
\end{proof}
\end{cor}

\begin{remark}
Theorem \ref{thm: Main} and Corollary \ref{cor: n_k = 1 for ret nest} give necessary conditions for the return times of a return
nest.  It is natural to ask what are sufficient conditions?  Moreover are these conditions different for polynomials with a
connected Julia set? with a disconnected Julia set? or for real polynomials? However, we do not address these questions at this
time, since it would require a more in depth and technical analysis of a tree with dynamics than we wish to present in this
paper.
\end{remark}

\subsection{Polynomials with Disconnected Julia Sets}\label{subs: Poly - Disc J}

We now turn to describing the standard puzzle of a polynomial.  First we recall some facts that we will use in both the connected
and disconnected cases.

Let $f$ be a polynomial of degree $d\geq 2$. Let $g$ denote Green's Function of $f$. We use $g$ to decompose the plane. Recall
that $g(z) \geq 0$ for all $z\in \bC$ and $g(z)= 0$ if and only if $z \in \Kjulia(f)$.  The functional equation $g(f) = d \cdot
g$ is satisfied by $f$ and $g$. The critical points of $g$ are the critical points of $f$ and the iterated pre-images of critical
points of $f$.  An \emph{equipotential} is a level set of $g$ of the form $\set{z \in \bC: \ g(z) = \lambda > 0}$. By the
functional equation, $f$ maps equipotentials to equipotentials. Equipotentials have a dynamical definition (see \cite{Branner}).
It follows that the tree with dynamics of a polynomial with a disconnected Julia set is a topological invariant.

Let $\varphi$ be some B\"ottcher function of $f$. If the Julia set of $f$ is connected, then $\varphi^{-1}: \bC \minus
\overline{\bD} \to \bC \minus \Kjulia(f)$ is a conformal isomorphism. An \emph{external ray} is a set of the form
$\varphi^{-1}\set{s e^{2 \pi i \theta}: 0 < s < \infty}$ for some $\theta \in \bR / \bZ$. An external ray is \emph{rational} if
$\theta \in \bQ / \bZ$.

We outline the dynamical decomposition of the plane using $g$ for a polynomial with a disconnected Julia set. Fix a polynomial
$f$ of degree $d \geq 2$ with a disconnected Julia set.  We distinguish all equipotentials whose grand orbit contains a critical
point of $f$. There are countably many such equipotentials, say $\set{E_l}_{l \in \bZ}$. Index them so that $g|E_l < g|E_{l-1}$,
$E_l$ is a Jordan curve for $l \leq 0 $, and $E_1$ is not a Jordan curve (so it contains a subset homeomorphic to a figure-8).
Let $H$ be the number of orbits of $\set{E_l}_{l \in \bZ}$ under $f$. If $f$ has $e$ distinct critical points that escape to
infinity, then $H \leq e$. It is possible that $H < e$ if $f$ has two escaping critical points $c$ and $c'$ such that $g(c) = d^n
g(c')$ for some $n \in \bZ$. From the functional equation and the indexing of $E_l$, it follows that $f(E_l) = E_{l-H}$ for all
$l$.  In fact, $H$ is equivalent to the $H$ that is used in trees with dynamics. Notice that $H$ is positive.

Define $ U_l = \set{z:  g(z) <  g|E_{l} }$. For $l \leq 0$, $U_l$ is a single \emph{topological disk} (defined as a simply
connected open set).  For all $l$, $U_l$ is the disjoint union of finitely many topological disks $\set{P_l^i}$. We define each
of these disks as a puzzle piece of $f$ at depth $l$, and $\puz_l = \set{P_l^i}$. The \emph{Branner-Hubbard puzzle} is $\puz =
(\puz_l)$ \cite[Ch.\ 3]{BH92}. It is known that the Markov properties hold for these puzzle pieces.

\begin{definition}\label{defn: TwD of Poly w/Disc}
Let $f$ be a polynomial of degree at least 2 with a disconnected Julia set. Let $\puz$ be the Branner-Hubbard puzzle of $f$. We
define the \emph{canonical tree with dynamics} of $f$  as $(\tree(\puz), \F_f)$.\\
\end{definition}

\begin{figure} [h]

\centerline{\includegraphics*{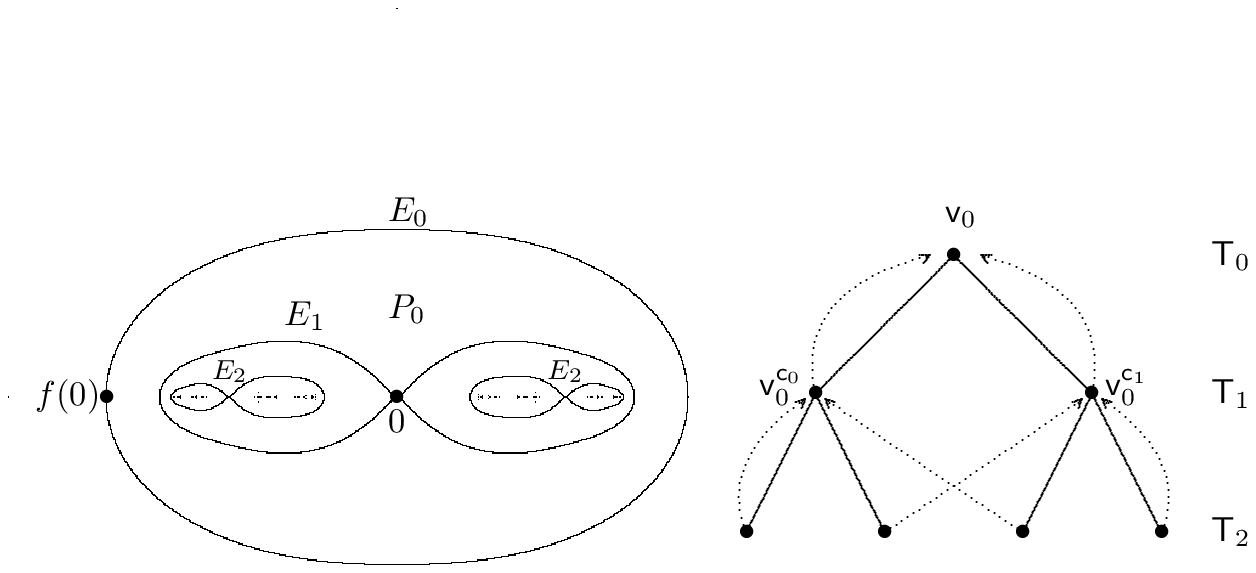}}

\caption{Equipotentials of a quadratic polynomial with a disconnected Julia set and its tree with dynamics.}

\label{fig: Quad Tree}
\end{figure}

\begin{example} \label{eg: disc quad TwD}
Let $f(z) = z^2 + c$ have a disconnected Julia set. There is exactly one escaping critical point of $f$, namely $0$. So $H = 1$.
The equipotential that contains $0$ is $E_1$. Each equipotential $E_{-l} = f^{l+1}(E_1)=\set{z: g(z)  = 2^{l+1} g(0)}$ is an
analytic Jordan curve for $l \in \bN$. Hence, $\tree_{-l}$ is a single vertex for $l \in \bN$. Since $0$ is a simple critical
point, $f$ is locally two-to-one near 0. Thus, $E_1$ is homeomorphic to a figure-8: two Jordan curves pasted at $0$. There are
exactly two components of $E_2 = f^{-1}(E_1) = \set{z: g(z) = 2^{-1}g(0)}$, one inside each loop of $E_1$. Thus, there are
exactly two puzzle pieces of $f$ at depth 1.  So in the tree with dynamics of $f$, $\tree_0 = \set{\mathsf{v}_0}$ and $\tree_1 =
\set{\mathsf{v}_0^{\child_0}, \mathsf{v}_0^{\child_1}}$. For $l >  1$, $f$ maps $E_l$ onto $E_{l-1}$.  This map is one-to-one,
since there are no critical points in the bounded components of $\bC \minus E_l$. Therefore each component of $E_l$ is
homeomorphic to a figure-8. It follows that there are exactly two components of $E_{l+1}$ nested inside each component of $E_l$;
one in each loop of the figure-8.  Thus if $\mathsf{v} \in \tree_l$ for some $l \in \bN$, then $\mathsf{v}$ has exactly two
children, which are mapped by $\F$ onto the children of $\F(\mathsf{v})$. It follows that $(\tree, \F)$ is the tree with dynamics
from Example \ref{eg: binary TwD}.

\end{example}

A \emph{child} of a puzzle piece $P\in \puz_l$ is a $P^{\child} \in \puz_{l+1} $ such that $P^{\child} \subset P $. Define  $A =
P \minus \bigcup \overline{P^{\child}}$ as an \emph{annulus of} $f$. So, there is a one-to-one correspondence between annuli and
puzzle pieces. Since $f$ preserves children, it follows that $f$ maps each annulus onto an annulus. Let $A_i = P_i  \minus
\bigcup \overline{P_i^{\child}}$ for $i=1,2$.  Then $f(A_1) = f(A_2)$ if and only if $f(P_1) = f(P_2)$.  Therefore, we could just
as well used annuli instead of puzzle pieces to define the puzzle.

\begin{remark}
Because of the correspondence between puzzle pieces and annuli, we could have defined the vertices of the tree as annuli.  Thus,
Definition \ref{defn: TwD of Poly w/Disc} gives the same tree with dynamics as \cite[Def.\ 3.7]{E03}.
\end{remark}

\subsection{Polynomials with Connected Julia Sets}\label{subs: Poly - Con J}

When $f$ is a polynomial with a connected Julia set, we choose a {Yoccoz puzzle} \cite{Hubbard_Loc_Con, Kiwi_Rat_Rays}. We
outline the construction here. All the equipotentials of $f$ are Jordan curves, so we use external rays to separate the Julia
set.

Fix a polynomial $f$ of degree $d \geq 2$ with a connected Julia set.  Choose a finite set $\alpha_1, \dots, \alpha_m \in
\julia(f)$ such that at least one rational ray lands on each $\alpha_i$, at least two distinct rational rays land on some
$\alpha_j$, and $f(\set{\alpha_i})\subset \set{\alpha_i}$. For quadratic polynomials, one usually chooses $\set{\alpha_i} =
\set{\alpha}$, where $\alpha$ is a repelling fixed point that is the landing point of at least two distinct external rays. In
general, we could choose a finite number of repelling periodic cycles. Choose $\lambda > 0$. The equipotentials $\set{z: g(z) =
d^l\lambda}$ and the pre-images of the external rays that land on $\set{\alpha_i}$ partition the plane. We define the Yoccoz
puzzle of $f$ determined by $\set{\alpha_i}$ and $\lambda$. Define $U_1 = \set{z: g(z) < \lambda } $, and $\Gamma_1$ as the
external rays that land on $\alpha_1, \dots, \alpha_m$ including the landing points. We define a sequence of open sets
$(U_l)_{l=1}^{\infty}$ and graphs $(\Gamma_l)_{l=1}^\infty$, by $U_{l+1} = f^{-1}(U_l)$ and $\Gamma_{l+1} = f^{-1}(\Gamma_l)$.
For technical reasons, define $U_{-l} = \set{z: g(z) < d^l\lambda } $ and $\Gamma_{-l} = \emptyset$ for $l \in \bN$. A connected
component $P_l^i$ of $U_l \minus \Gamma_l$ is a \emph{puzzle piece of $f$ at depth} $l$, and $\puz_l = \set{P_l^i}$. Define the
Yoccoz puzzle by $\puz = (\puz_l)$.  Note that what we call depth $l$, previous authors call depth $l-1$. The Markov properties
hold for pieces at any level $l \geq 1$ \cite[\S12]{Kiwi_Rat_Rays}.

Puzzle pieces at depth $l \leq 0$ have not been defined in the existing literature. Each of these puzzle pieces is bounded by an
equipotential, but not by any external rays. Hence, there is only one puzzle piece at any depth $l \leq 0$, which is a
topological disk. It is easy to check the Markov properties hold for these pieces.  A puzzle piece at depth 1 may have
exceptional dynamics. Let $\puz_0 = \set{P_0}$. If $P$ is a puzzle piece at depth $1$, then $f(P) \subset P_0$. However in
general, $f(P) \neq P_0$. Even if $f(P) = P_0$, the map $f:{P} \to {P}_0$ is not necessarily proper. From Lemma \ref{lem: TwD of
Puz}, we have $F_f(P) = P_0$. The exceptional nature of these pieces does not affect return times.

\begin{definition}\label{defn: TwD of Poly w/Con}
Let $f$ be a polynomial with a connected Julia set. Let $\puz$ be some Yoccoz puzzle of $f$. We define the \emph{standard tree
with dynamics with respect to} $\puz$ as $(\tree(\puz), \F_f)$.
\end{definition}

\begin{remark}
The tree with dynamics of a polynomial with a connected Julia set depends on the choice of $\puz$.  Specifically on the
$\set{\alpha_i}$ that are chosen in the construction $\puz$. Hence the tree with dynamics of a polynomial with a connected Julia
set is not canonical.  However it does not depend on the choice of $\lambda$.
\end{remark}

\begin{figure}[hbt]
\begin{center}

\includegraphics{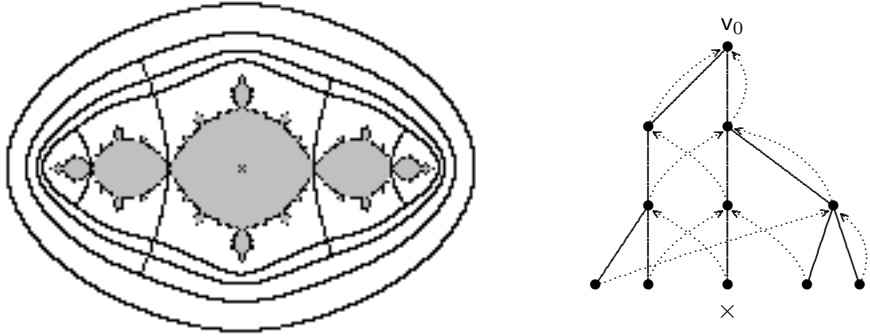}

\caption{A Yoccoz puzzle of the Julia set of $z^2 - 1$ and its tree with dynamics. The end which corresponds to the nest of the
critical point $(0)$ is indicated by the symbol ${\times}$.}

\label{fig: TwD of Yoccoz Puzzle}

\end{center}
\end{figure}

\begin{example}
Let $f(z) = z^2 +c$ have a connected Julia set. Assume exactly 2 external rays land on some fixed point $\alpha$.  We outline the
construction of the first few levels of a Yoccoz puzzle of $f$ and the associated tree with dynamics, see Fig.\ \ref{fig: TwD of
Yoccoz Puzzle}. We choose $\set{\alpha_i} = \set{\alpha} $, and any $\lambda >0$. By definition, $U_{-l} = \set{z: g(z) < 2^{l+1}
\lambda}$ and $\Gamma_{-l} = \emptyset$ for $l  \in \bN$.  Thus $\puz_{-l}$ is a single puzzle piece $P_{-l}$, and $\tree_{-l} =
\set{\mathsf{v}_{-l}}$ for $l  \in \bN$.  We have $\Gamma_1 =\set{\frac{1}{3}, \frac{2}{3}}$, and these two external rays
separate $U_1$ into two puzzle pieces, say $P_1^{0}$ which contains 0, and $P_1^1$. We have $f(P_1^1) \subsetneq P_0$ and
$f(P_1^0) = P_0$.  However, the map $f: P_1^0 \to P_0$ is not proper. Next, $\Gamma_2 =\Gamma_1 \cup \set{\frac{1}{6},
\frac{5}{6}}$. So $P^0_1$ has two children, say $P^0_2 \ni 0$ and $P_2^1$, and $P^1_1$ has only one child $P_2^2$. So $f(P_2^0) =
P_1^1$, $f(P_2^1) = P_1^0$, $f(P_2^2) = P_1^0$.   Finally, $\Gamma_3 =\Gamma_2 \cup \set{\frac{1}{12}, \frac{11}{12},
\frac{5}{12}, \frac{7}{12}}$. This gives $P_2^1$ and $P^2_2$ two children each, and $P_2^0$ one child.
\end{example}

In the connected case, an \emph{annulus of $f$ at depth $l$} is the difference of nested puzzle pieces $P^{\parent} \minus
\overline{P}$ where $P$ is a puzzle piece at depth $l$. A complication can arise. An annulus $A = P^{\parent} \minus
\overline{P}$ is called \emph{degenerate} if $\bdr P\cap \bdr P^{\parent} \neq \emptyset$. A degenerate annulus is not doubly
connected, but is the union of topological disks.  This complication is a serious concern in modulus estimates (see
\cite[\S9]{Hubbard_Loc_Con}), but does not affect return times.  Annuli are mapped onto annuli by $f$, except possibly for those
at depth $1$. Let $A_i = P_i^{\parent} \minus \overline{P_i}$ for $i=1,2$. If $f(A_1)  = A_2$, then $f(P_1^{\parent}) =
P_2^{\parent}$ by the disjointness property of puzzle pieces.

\begin{remark}
As in the disconnected case, if we use annuli instead of puzzle pieces, we obtain the same tree with dynamics of $f$.
\end{remark}

\subsection{Properties of Polynomial Trees with Dynamics}\label{subs: Poly - Poly TwD}

We note some properties satisfied by any tree with dynamics of a standard puzzle of a polynomial. The following properties follow
immediately from the construction of a standard tree with dynamics of a polynomial (Definitions \ref{defn: TwD of Poly w/Disc}
and \ref{defn: TwD of Poly w/Con}). Compare to \cite[Def. 3.2]{E03}.

\begin{prop}\label{prop: Poly TwD properties}
If $(\tree, \F)$ is a tree with dynamics of a polynomial with respect to a standard puzzle, then the following hold.
\begin{enumerate}

    \item Every vertex of $\tree$ has at least one child.  That is, $\tree$ has no leaves.

    \item Every vertex of $\tree$ has only finitely many children. That is, $\tree$ is locally finite.

    \item There is a root of $\tree$.

    \item The integer $H$ such that $\F(\tree_l ) = \tree_{l-H}$ is positive. Moreover $H = 1$ for any a tree with dynamics
        of a polynomial with a connected Julia set.

\end{enumerate}

\end{prop}

In this paper, the construction of a tree with dynamics from a puzzle is purely set theoretic.  We only use the Markov properties
of the puzzle. We never used the topological properties of puzzle pieces of a polynomial.  Thus, the above proposition is a list
of the properties of a tree with dynamics of a polynomial, which follow from the set-theoretic properties of a standard puzzle.
These standard puzzles have some useful additional properties; all of their puzzle pieces are open, simply connected, and
pre-compact. Also, a polynomial restricted to one of its puzzle pieces is a proper map, with at most finitely many exceptions.
Notice that we never use these additional properties in this paper.

Axioms for the tree with dynamics of a polynomial with a disconnected Julia set are known. These axioms are complete in the sense
that they are necessary and sufficient conditions for a tree with dynamics to be the tree with dynamics of a polynomial with a
disconnected Julia set.  These axioms were first given by the author in \cite[Def.\ 4.7]{E03}, where their necessity was
demonstrated. They were also shown to be sufficient for certain trees with dynamics \cite[Thm.\ 9.4]{E03}. L.~DeMarco and
C.~McMullen showed that these axioms are sufficient (using different notation) \cite{DM-MC-tdp}.  To define these complete axioms
we need to keep track of some topological information.  Specifically, the local degree of the polynomial restricted to each of
its puzzle pieces. A formula for the number of children of each puzzle piece can then be deduced from the Riemann-Hurwitz formula
for domains \cite[Def.\ 4.6]{E03}.

\begin{remark}
To the best of the author's knowledge complete axioms for a tree with dynamics of a polynomial with a connected Julia set are not
known.  There are two obvious differences from the disconnected case. The number of children of the root $\mathsf{v}_0$ can be
arbitrarily high (depending on the choice of $\alpha_1, \dots, \alpha_m$).  The axioms for vertices at level 1 are not clear,
since the polynomial restricted to a puzzle piece at level 1 is not necessarily proper. For all other levels, the polynomial
restricted to each puzzle piece is proper.  Thus, the conditions of \cite[Def.\ 4.6]{E03} must hold.
\end{remark}

We now describe how the dynamics of a nest corresponds to the dynamics on a Julia set. Let $f$ be a polynomial and $\puz$ a
puzzle of $f$.  For $z_0 \in \Kjulia(f)$, a \emph{nest of} $z_0$ is a nest of puzzle pieces $(P_l)$ such that $z_0 \in
\overline{P_l}$ for every $l$. If we regard a return of $z_0$ to $P_{l+1}$ as ``closer'' than a return to $P_l$, then the return
times of a return nest can be thought of as generalized closest return times.

Let $(P_l)$ be a nest of puzzle pieces from the Branner-Hubbard puzzle of $f$ (Definition \ref{defn: TwD of Poly w/Disc}).  Then
$\bigcap_{l \in \bN} \overline{P_l}$ is a connected component of $\Kjulia(f)$.  So each point $z_0 \in \Kjulia(f)$ lies in a
unique nest. Thus the dynamics of $z_0$ are coded by the dynamics of a unique end of the tree with dynamics.  Moreover, if
$\julia(f)$ is a Cantor set, then there is a one-to-one correspondence between points of $\julia(f)$ and ends of the tree with
dynamics.

Let $\puz$ be a Yoccoz puzzle of $f$ (Definition \ref{defn: TwD of Poly w/Con}).  Let $z_0 \in \Kjulia(f)$.  If $z_0$ is not a
pre-image of one of the $\alpha_i$ used to define the Yoccoz puzzle, then $z_0$ lies in a unique nest $(P_l)$. Although there is
no obvious geometric interpretation of $\bigcap_{l \in \bN} \overline{P_l}$. If $z_0$ is the pre-image of one of the $\alpha_i$
used to define the Yoccoz puzzle, then $z_0 \in \bigcap_{l \in \bN} \overline{P_l}$ for several nests $(P_l)$.  Thus the dynamics
of every point of $\Kjulia(f)$ are coded by an end of the tree with dynamics.  For all but countably many of points of
$\Kjulia(f)$, this coding is unique.

A nest of puzzle pieces that contains a critical point is of special interest.

\begin{definition}\label{defn: principal nest}
Let $(P_l)$ be an extended nest of puzzle pieces of a polynomial $f$. Suppose that there is a critical point of $f$ in every
$P_l$. A \emph{principal nest} is a minimal return nest of $(P_l)$.
\end{definition}

In terms of Markov properties, there is no difference between principal nests and return nests.  The return times of both are
$r(k)$-bonacci numbers. However, topologically they are different. A principal nest can be to used estimate the moduli of the
annuli of its nest \cite[Ch.\ 4.5]{BH92}. A general return nest cannot be used in this way.

The concept of a principal nest appeared in the literature before the terminology was set. Branner and Hubbard introduced
tableaux \cite[Ch.\ 4.2]{BH92}.  For a polynomial with exactly one persistent critical point, a tableau keeps track of the return
times of a principal nest of the persistent critical point with $l(0) =0 $. Later the zig-zag pattern that defined the returns
was called the \emph{critical staircase} of the tableau \cite{Branner}.  Lyubich suggested a choice for $l(0)$, and we have
generalized the terminology he introduced in \cite[\S3]{Lyubich-Quadratics}.  A standard puzzle of polynomial with a unique
persistent critical point has a principal nest, which depends only on the choice of $l(0)$ (Proposition \ref{prop: form min ret
chain}).

A principal nest can sometimes be used to find domains where a polynomial is renormalizable in the sense of Douady and Hubbard
\cite{DH-Poly_like}. Let $(P_{l(k)})_{k \in \bZ}$ be a principal nest of a polynomial $f$.  Then $f|P_{l(k)}$ is not univalent
for any $k$.  For a Branner-Hubbard puzzle, it is clear that $(f^{n_k}; P_{l(k)}, P_{l(k-1)})$ is polynomial-like (of some degree
$\geq 2$)  for every $k \in \bZ$.  A polynomial of the form $f(z) = z^d+c$ ($d\geq2$, $c\in \bC$) is called \emph{unicritical}.
Lyubich showed that for a Yoccoz puzzle of a unicritical polynomial, there is a choice of $l(0)$ for the principal nest such that
$(f^{n_k}; P_{l(k)}, P_{l(k-1)})$ is polynomial-like for every $k \geq 1$ \cite[Prop.\ 3.1]{Lyubich-Quadratics}.

\vspace{.2in}

\textit{Acknowledgements.} Portions of this paper are based upon my Ph.\ D.\ thesis. I am indebted to Ricardo P{\'e}rez-Marco for
his guidance.

The computer generated images of Julia sets were created using \emph{Mandel} by Wolf Jung
(\href{http://www.mndynamics.com}{http://www.mndynamics.com}), and \emph{DH-Drawer} by Arnaud Cheritat
(\href{http://www.math.univ-toulouse.fr/~cheritat/}{http://www.math.univ-toulouse.fr/{$\sim$}cheritat/}).

\providecommand{\bysame}{\leavevmode\hbox to3em{\hrulefill}\thinspace}

\end{document}